\documentclass{article}

 \usepackage{amsfonts}

\newtheorem{theo}{Theorem}[section]

\newtheorem{defin}[theo]{Definition}

\newcommand{\Reali}{\ensuremath{\mathbb{R}}}

\begin{document}
\title{An Approximation to the Gr\"obner Basis of Ideals of Perturbed Points: Part I.}
\author{Claudia Fassino\thanks{Dipartimento di Matematica, Universit\`a di Genova, via Dodecaneso 35, 16146 Genova, Italy
 ({\tt fassino@dima.unige.it}).}}
\date{}
\maketitle
\begin{abstract}
We develop a method for approximating the Gr\"obner basis of the
ideal of polynomials which vanish at a finite set of points, when
the coordinates of the points are known with only limited
precision. The method consists of  a preprocessing phase  of the
input points to mitigate the effects of the input data
uncertainty,  and of a new ``numerical" version of the
Buchberger-M\"oller algorithm to compute an approximation
$\overline{GB}$ to the exact Gr\"obner basis. This second part is
based on a threshold-dependent procedure for analyzing from a
numerical point of view the membership of a perturbed vector to a
perturbed subspace. With a suitable choice of the threshold, the
set $\overline{GB}$ turns out to be a good approximation to a
``possible" exact Gr\"obner basis or to a basis which is an
``attractor" of the exact one. In addition, the polynomials of
$\overline{GB}$ are ``sufficiently near"  to the polynomials of
the extended basis, introduced by Stetter, but they present the
advantage that $LT(\overline{GB})$ coincides with the leading
terms of a ``possible" exact case. The set of the preprocessed
points, approximation to the unknown exact points, is a pseudozero
set for the polynomials of $\overline{GB}$.
 \end{abstract}
\section{Introduction}

Let $P=K[x_1,\dots,x_s]$ be the polynomial ring in $s$
indeterminates over a field $K$. In this paper we develop a method
for approximating the Gr\"obner basis of the ideal ${\cal I}
\subset P$ of polynomials which vanish at a finite set of points,
 when the coordinates of the points are known with only
limited precision. In particular we analyze the case $K=\Reali$.

In the exact case, that is when we deal with a set of unperturbed
points $\cal P$, the problem  has been deeply analyzed by several
authors (see for example \cite{AKRobbiano05}, \cite{KrRob00},
\cite{KrRob05}, \cite{MMM95}) and the Gr\"obner basis of the ideal
$\cal I({\cal P})$ of the points $\cal P$ can be computed, for
example, by the Buchberger-M\"oller algorithm (\cite{ABR00},
\cite{BM},  \cite{KrRob05}).  It is well known that, given a term
ordering, each set of $m$ points $\cal P$ corresponds to a closed
set ${\cal N}=\{t_1,\dots,t_m\}$ such that each polynomial $g$ of
the Gr\"obner basis $GB$ of $\cal I({\cal P})$ has the form
$$g=t-\sum_{t_i < t} c_i t_i,$$
with leading term $t$ and suitable coefficients $c_i$. The set
$\cal N$, called {\bf normal set}, is the basis of the quotient
ring $P/\cal I({\cal P})$ and the leading term $t$  belongs to the
corner set ${\cal C}[{\cal N}]$.

Nevertheless, it is also well known that the problem of computing
the Gr\"obner basis of an ideal of points is ``discontinuous"
(\cite{KrRob05}, \cite{Ste04}), i.e. small perturbations of the
points can cause structural changes in the basis, as illustrated
in Example 1.1. For this reason exact methods applied blindly to
perturbed ``real-world" data can produce meaningless results.

\vspace{3 mm}\noindent {\bf Example 1.1.} Given the term ordering
DegLex and the exact points   $p_1=(1,1)$, $p_2=(3,2)$ and
$p_3=(5,3)$, the Gr\"obner basis of the ideal $\cal I$ of
polynomials in $\Reali[x,y]$ vanishing at these points is
\begin{eqnarray*}
GB= \left\{
\begin{array}{l}
x-2y+1, \\
 y^3-6y^2+11y -6,
\end{array}
\right. \;\;{\rm with \; normal\; set \;} {\cal N}
=\{1,\,y,\,y^2\}.
\end{eqnarray*}
If  data errors  perturb the point $p_3$ so that the set of input
points is  $\widehat{\cal P} =\{p_1,\; p_2, \; (5.1,3)\}$, then
the Gr\"obner basis $\widehat{GB}$ of the ideal ${\cal
I}(\widehat{\cal P})$ is:
\begin{eqnarray*}
\widehat{GB}=
\left\{
\begin{array}{l}
y^2-20x+37y-18, \\
 xy-43x+81 y -39,\\
 x^2-90.1x+172.2y-83.1,
\end{array}
\right.\;\;{\rm with \; normal\; set \;} \widehat{\cal N}
=\{1,\,y,\,x\}.
\end{eqnarray*}
This is an obvious result since the points $\widehat{\cal P}$ are
not aligned and so  $\widehat{GB}$ structurally  differs from
$GB$. Nevertheless, since the points $\widehat{\cal P}$ are not aligned
because of data errors, from the computational point of
view it is more interesting to approximate the Gr\"obner basis
$GB$ than to compute $\widehat{GB}$, given as input the set
$\widehat{\cal P}$. \hfill{$\diamondsuit$}

\vspace{3 mm} In the general case, if the input points arise from
real-world data, their coordinates are approximations to the
unknown exact values, and only an error estimation is known. For
this reason, each point $\widehat p$ of $\widehat{\cal P}$
represents a ``cloud" of points. More precisely, given an input
point $\widehat p$  known with uncertainty $s_0$, each point
$\widetilde p$ which differs from $\widehat p$ componentwise by
less than $s_0$ can be chosen as an input point computationally
equivalent to $\widehat p$, and so it will be called ``{\bf
admissible input point}". Analogously, given a set of perturbed
input points $\widehat {\cal P}$, an ``{\bf admissible input set
}" is a set consisting of admissible input points computationally
equivalent to $\widehat{\cal P}$. Obviously, given an error
estimation $s_0$, there exist several admissible input sets. These
remarks combined with the structural  discontinuity show that the
computation of the Gr\"obner basis of ideals of perturbed points
is a very tricky problem. In fact, since the exact points are
unknown and the Gr\"obner basis can change choosing different
admissible input sets, it can be difficult to identify the case to
be approximated. For this reason, we define some  criteria to
choose the ``reference" Gr\"obner basis to approximate, pointing
out its structure, that is its leading terms and the supports of
its polynomials.

We define  ``{\bf reference basis}"  a Gr\"obner basis $GB$
generating an ideal of polynomials vanishing at some admissible
input set, whose structure represents all the admissible input
sets. Since  the structure of a  Gr\"obner basis and the elements
of its normal set are closely connected, we require that the
normal set of a reference basis can be associated to all the
admissible input sets, that is, analogously to the exact case, we
require that its terms provide independent vectors if evaluated on
the points of any admissible input set. A normal set with this
property will be called ``{\bf reference normal set}". Note that a
reference basis could coincide with the exact basis, since the set
of exact points belongs to the cloud of the admissible input sets.

In order to approximate a reference  basis, our method computes a
set of polynomials $\overline{GB}$, whose supports belong to
$\overline{\cal N} \cup {\cal C}[\overline{\cal N}]$, where
$\overline{\cal N}$ is a reference normal set. In some particular
case, $\overline{GB}$ is a reference Gr\"obner basis, and, in
general, it is  a good approximation to a reference basis. In
fact, we show that, in the general case,
$\overline{GB}=\{\overline g_1,\dots,\overline g_n\}$ {\bf
structurally corresponds} to a reference basis
$GB=\{g_1,\dots,g_n\}$, that is there is an one-to-one
correspondence between $\overline{GB}$ and $GB$ which preserves
the leading terms and the supports of each polynomial. In
addition, for each $i=1\dots n$, $\overline g_i \in\overline{GB}$
and $g_i \in GB$ have ``similar" coefficients and we estimate the
difference between $\overline g_i$ and $g_i$. Moreover, we show
that the elements of $\overline{GB}$ are ``near" to the
polynomials of the extended basis $GB_E$ introduced by Stetter in
\cite{Ste97}.

Unfortunately,  some ``{\bf degenerate}" cases can occur, when
there are no reference bases, that is when all the Gr\"obner
bases, with reference normal set, are associated to points which
are not admissible input points. Nevertheless, in these degenerate
cases the numerical tests point out that also the normal set of
the exact points cannot represent, from a numerical point of view,
all the admissible input sets. In addition, the tests suggest that
the computed normal set $\overline{\cal N}$ seems to be the
numerically stable normal set ``nearest" to the exact one. The
analyzed degenerate examples are not presented in this paper since
the degenerate cases will be studied in a future work, whereas in
this article we shall limit ourselves to ``{\bf non-degenerate}"
cases.

Our method consists of two parts. First of all, the input points
are preprocessed, in order to mitigate some negative effects of
the data errors. The possible ``splittings" of the coordinates due
to the errors are removed and the preprocessing algorithm computes
a new admissible input set $\overline{\cal P}$, by replacing with
a unique suitable value the perturbed coordinates which differ by
less than the error estimation. The preprocessed points turn out
to be a set of pseudozeros \cite{Ste199} for the polynomials of
$\overline{GB}$.

The second part of our algorithm constitutes the main innovative
part of this paper. It is a numerical version of the
Buchberger-M\"oller algorithm which  computes a set of polynomials
$\overline {GB}$, approximation to a reference basis. The
numerical Buchberger-M\"oller algorithm is based on the following
notion of ``approximate" linear dependence of vectors. From the
numerical point of view, a vector $v$ depends on the vectors
$\{v_1,\dots,v_k\}$ if it is ``sufficiently" near to an element of
the subspace ${\rm Span}\{v_1,\dots,v_k\}$. Since this
``numerical" dependence of vectors is defined choosing a threshold
$\epsilon$, the numerical version of the Buchberger-M\"oller
algorithm is threshold-dependent. We have to choose a suitable
value of the threshold for computing a good approximation to a
reference  basis. The threshold choice is one of the trickier
aspects of the implementation of our method.

\vspace{3 mm} This paper is organized as follows. In Section 2 we
describe the threshold-dependent approximation algorithm and in
Section 3 we present an heuristic approach for the choice of this
threshold, whereas a more detailed  analysis is left to future
works. The approximation properties of the computed set
$\overline{GB}$ are analyzed in Section 4, and a comparison with
some results presented by Stetter is described in  Section 5. In
Section 6 we report two examples in order to illustrate and
clarify explicitly the behaviour of our method. Appendix A
contains a description of the strategy for preprocessing the input
points used in our examples and Appendix B presents two error
upper bounds useful in the estimation of  the goodness of the
basis approximation.

\vspace{3 mm}\noindent {\bf Notation.} Later on ``G-basis"  means
``Gr\"obner basis". We denote the points with the letter $p$ and
the tuple of points with  ${\cal P}=\{p_1,p_2\dots\}$. If
$g(x_1,\dots,x_s) $ is a polynomial in $\Reali[x_1,\dots,x_s]$,
$g(p)$ is the evaluation of $g$ at  $p$  and $g(\cal{P})$ is the
vector whose components are the values of $g$ on the elements of
$\cal{P}$. In addition we denote with $\|v\|_2$ the 2-norm of a
vector $v$.

\section{The approximation method}
In this section we describe our method for approximating a
reference basis with a set of polynomials $\overline{GB}$, when we
deal with perturbed data in non-degenerate cases, that is when a
reference basis exists.

Our method computes  the set of polynomials $\overline{GB}$,
working on  a set of points obtained after preprocessing the
perturbed input data. The goal of the preprocessing phase is to
minimize the damage caused by possible ``splittings" of the data
due to the input perturbations, replacing by a unique value the
coordinates which differ by less than the error estimation, since
they are coincident from the computational point of view. The
preprocessed set $\overline{\cal P}$ depends on the preprocessing
technique; the strategy used in our tests is described in Appendix
A.

Even if the exact  Buchberger-M\"oller (BM) algorithm \cite{BM}
applied to the preprocessed data sometimes computes a good
approximation to a reference basis, this is not true in general.
The simple preprocessing of the input data is not sufficient to
obtain a numerically  significant set of polynomials, and so it is
necessary to develop a new approximation algorithm.

It is well known that, given a term ordering $\sigma$ and a set of
$m$ points ${\cal P}=\{p_1,\dots,p_m\}$, not perturbed by errors,
the BM algorithm computes the $\sigma$-Gr\"obner basis $GB$
 of the ideal of the polynomials which vanish at the set
$\cal P$. At each step, after choosing the term $t$, the BM
algorithm checks if the vector $t({\cal P})$ is linearly
independent of the vectors $\{t_1( {\cal P}),\dots,t_k({\cal
P})\}$, where $\{t_1,\dots,t_k\}$ is the subset of the normal set
computed at the previous steps, and $t>_\sigma t_i$,
$i=1,\dots,k$. If this is the case, $t$ is a new term of the
normal set. Otherwise, the BM algorithm builds a new polynomial
$g$ of the G-basis $GB$ such that $g=t-\sum_{i=1}^k c_i t_i$. The
coefficients $c_i$, $i=1,\dots,k$, are the components of the
vector $c$ that satisfies
\begin{eqnarray}\label{sistema_esatto}
 M_k({\cal P}) c = t({\cal P}),
\end{eqnarray}
where $M_k({\cal P})$ is the matrix $ [t_1({\cal
P}),\dots,t_k({\cal P})]$, whose columns are linearly independent,
by construction.

The above test of linear dependence is crucially affected by even
very small variations of the input data. Therefore, when we deal
with approximate input points, different choices of admissible
data may lead to different choices in the BM algorithm.
Nevertheless, if the data are affected by errors, it is possible
to check the linear dependence of perturbed vectors from a
numerical point of view. Intuitively, a perturbed vector $v$
belongs, in some approximate sense, to a subspace $W$, spanned by
perturbed vectors, if $v$ is ``sufficiently near" to an element of
$W$. This idea of ``numerical" membership of a vector to a
subspace, formalized and analyzed by the author in \cite{Fas05},
is based on the following definition.
\begin{defin}\label{app}
Given a subspace $W$ and  a threshold $\epsilon>0$, a vector $v$
``numerically" lies in $W$ if
\begin{eqnarray*}
 \frac{\|r \|_2}{\|v\|_2}\le \epsilon,
\end{eqnarray*}
where $r$ is the component of $v$ orthogonal to $W$. In this case
we write $v\in_\epsilon W$.
\end{defin}

By exploiting this concept of ``numerical" membership of a
subspace, we develop the numerical {\bf A}pproximate {\bf
B}uchberger {\bf M}\"oller algorithm (denoted by ABM) for
computing the approximation $\overline{GB}$ to a reference basis.
Working on the preprocessed input points $\overline{\cal P}$, the
ABM algorithm checks, at each step, if a perturbed vector lies
numerically in a subspace spanned by a perturbed basis. If the
vector ``numerically" lies in the subspace, then the algorithm
builds a new polynomial of $\overline{GB}$, even if  the vector
does not exactly belong to the subspace. Otherwise the ABM
algorithm inserts a new element into the normal set
$\overline{\cal N}$.

The ABM algorithm can be described as follows. We check the
 numerical membership of a vector $v$ to a subspace $W={\rm Span}\{w_1,\dots,w_k\}$
computing the residual $r \in W^\perp$ of the least square problem
$[w_1\dots w_k] x = v$, since, in this way,  we can more easily
estimate the difference between $\overline {GB}$ and a reference
basis, as shown in Section 4.

\newpage
\noindent{\bf The ABM Algorithm.}

\vspace{2 mm}\noindent{\bf Input.} The term ordering $\sigma$, the
preprocessed input points $\overline{\cal P}=\{\overline
p_1,\dots,\overline p_m\}$ and the threshold $\epsilon$.

\vspace{2 mm}\noindent{\bf Output.} The  set of polynomials
$\overline{GB}$  and the set of power products $\overline{\cal
N}$.

\vspace{2 mm}\noindent Let  $\overline t_1,\dots,\overline t_k $
be the terms of $\overline{\cal N}$ computed at the previous steps
and let $\overline t$ be the current power product, $\overline
t>_\sigma \overline t_1,\dots,\overline t_k $, chosen analogously
to the BM algorithm.
\begin{enumerate}%
 \item Given the matrix  $\overline M_k(\overline{\cal P}) =[\overline t_1(\overline{\cal P}), \dots,\overline t_k(\overline{\cal P})]$, solve the least square problem  \begin{eqnarray}\label{pbminq}
\overline M_k(\overline{\cal P})\overline c = \overline
t(\overline{\cal P}).
\end{eqnarray}%
 \item Let $\overline r=\overline t(\overline{\cal P}) - \overline M_k \overline c $ be the residual of the problem (\ref{pbminq}),  then%
 \begin{itemize}%
 \item if  $ \|\overline r \|_2 / \|\overline t(\overline{\cal P})\|_2> \epsilon, $
put the term $\overline t$ into the set $\overline {\cal N}$;
 \item otherwise, put the polynomial  $\overline g = \overline t - \sum_{i=1}^k
\overline c_i  \overline t_i$ into $\overline{GB}$.\hfill $\diamondsuit$%
 \end{itemize} %
\end{enumerate}

\vspace{3 mm} Note that, analogously to the exact case, the
supports of the polynomials of $\overline{GB}$ consist of terms of
$\overline {\cal N} \cup {\cal C}[\overline{\cal N}]$ and the
leading terms  belong to ${\cal C}[\overline{\cal N}]$.

Obviously, the output sets $\overline {GB}$  and $\overline{\cal
N}$ of the ABM algorithm depend on the value of $\epsilon$ and we
have to choose a suitable value of the threshold so that
$\overline{GB}$ is a good approximation to a reference basis. A
detailed analysis of the choice of suitable thresholds will be
presented in a future work. In the next section we illustrate an
heuristic approach based on the analysis of  some properties of
$\overline{\cal N}$ and $\overline{GB}$.

\section{The threshold choice}
Intuitively, if we choose too small a value of $\epsilon$, the ABM
algorithm might not recognize numerically dependent vectors and so
the set $\overline {GB}$ could be very sensitive to data
perturbations. Vice versa, if the value of the threshold is too
large a vector may be considered as numerically belonging to a
subspace even if it is too ``far" from it. In this section we
present an heuristic strategy for choosing a suitable value of the
threshold analyzing the computed output sets. If the value of
$\epsilon$ does not produce satisfactory output sets, we repeat
the procedure with a new threshold.

As mentioned earlier, $\overline{GB}$ is a good approximation to a
reference basis $GB$ if each polynomial of $\overline{GB}$
corresponds to a polynomial of $GB$ with the same support and
similar coefficients and vice versa. In order to obtain the same
structure of a reference basis, since we have excluded degenerate
cases, we require that $\overline{\cal N}$ is the normal set of a
reference basis, and so we have to choose a value of $\epsilon$
such that the following properties hold.
\begin{description}
 \item{P1.} We require that  $\overline{\cal N}$ has  $m$ elements.
In this case, since the ABM algorithm analyzes the power products
using the same strategy of the exact BM algorithm, $\overline
{\cal N}=\{\overline t_1,\dots, \overline t_m\}$ turns out to be,
by construction, a closed set. For this reason, as is well known
\cite{Ste04}, $\overline{\cal N}$ is the basis of the quotient
ring $P/{\cal I}(\widetilde {\cal P})$ for each polynomial ideal
${\cal I}(\widetilde {\cal P})$ with $m$ simple zeros $\widetilde
{\cal P}$ such that $\{\overline t_1(\widetilde {\cal P}), \dots
\overline t_m(\widetilde {\cal P})\}$ are linearly independent
vectors, that is $\overline{\cal N}$ is a normal set associated to
such set  $\widetilde{\cal P}$.
 \item{P2.} We require that for each term $\overline t$   the
relative residual $\|\overline r\|_2/\|\overline t(\overline{\cal
P})\|_2$ of the problem (\ref{pbminq}) with right side $\overline
t(\overline{\cal P})$  must be sufficiently greater or less than
the threshold, that is well separated from it. This condition
guarantees that  small data perturbations do not affect each
residual sufficiently to change its position with respect to
$\epsilon$ and so $\overline{\cal N}$ turns out to be invariant
for admissible input sets. It follows that the elements of
$\overline{\cal N}$ provide independent vectors if evaluated on
admissible input sets, and so, if P1 holds, $\overline{\cal N}$ is
a normal set for the admissible input sets, i.e. it is a reference
normal set.
 \item{P3.} We require that, given the term ordering $\sigma$,
$\overline{\cal N}$ is the normal set of a $\sigma$-Gr\"obner
basis $GB$ of an ideal of admissible input points, that is
$\overline {\cal N}$ is the normal set of the reference G-basis
$GB$.
\end{description}
The properties P1 and P2 can be directly checked on the set
$\overline{\cal N}$ and on the computed relative residuals,
whereas it can be more difficult to investigate if $\overline{\cal
N}$ satisfies P3, since several cases can occur.

If the  following system, consisting of the polynomials of
$\overline{GB}$,
\begin{eqnarray}\label{sistema}
\left \{ \begin{array}{rcl}
\overline g_1 &=& 0 \\
\vdots & & \\
\overline g_n &=& 0,
\end{array}\right.
\end{eqnarray}
vanishes at $\overline{\cal P}$ then  $\overline {GB}$ is the
$\sigma$-Gr\"obner basis of the ideal ${\cal I}(\overline{\cal
P})$ with reference normal set $\overline{\cal N}$, since in this
case the ABM and BM algorithms have the same behaviour. Since
$\overline{\cal P}$ is an admissible input set, then
$\overline{GB}$ is a reference basis and  the property P3 is
obviously satisfied.

Otherwise we can check the property P3 using the following
intuitive idea. The value of the threshold is suitably chosen if
$\overline {GB}$ approximates the  reference basis
$GB=\{g_1,\dots,g_n\}$, with normal set $\overline{\cal N}$,
corresponding to a set of points $\cal P$ near to $\overline {\cal
P}$, that is if the function:
\begin{eqnarray}\label{funzione}
 \overline F: p  \mapsto (\overline g_1(p),\dots,\overline g_n(p)),
\end{eqnarray}
is a good  approximation to the function:
\begin{eqnarray*}
 F: p  \mapsto ( g_1(p),\dots,g_n(p)), \end{eqnarray*}
which vanishes at $\cal P$.
 For this reason, if the threshold
is suitably chosen,  we can suppose that  Newton method \cite{DBA}
works analogously if applied to the polynomial system $\overline
F(p)=0$ or to $F(p)=0$. Since Newton iteration applied to the
equation $\overline F(p)=0$ computes descent directions for
$\|\overline F(p)\|_2$ we can assume that it computes descent
directions for $\| F(p)\|_2$ too. For this reason the set
$\widetilde{\cal P}^+$, computed by Newton iteration working on
$\overline F$ with initial points $\overline {\cal P}$, can be
considered a good approximation to $\cal P$. It follows that the
property P3 is ``numerically" verified if the set $\widetilde{\cal
P}^+$ is sufficiently ``near" to $\overline{\cal P}$.

Since the properties P1--P3 must be checked on the computed sets
$\overline{\cal N}$ and $\overline{GB}$, we need a value of the
threshold to begin the computation. The first value of the
threshold can be chosen as follows. We apply the exact BM
algorithm to the preprocessed points $\overline{\cal P}$ and we
compute for each term $t$ of the exact normal set $\overline{\cal
N}_0$ the relative residual of the least square problem with right
side $t(\overline{\cal P})$. We choose the first value of
$\epsilon$ in order to eliminate from $\overline{\cal N}_0$, by
means the ABM algorithm, the power products with too small
relative residuals.

The following example illustrates the effects of different thresholds.

\vspace{3 mm} \noindent{\bf Example 3.1}

Let $\widehat{\cal P}=\{(1,\;1),\;\;(3,\;2),\;\;(5.1,\;3)\}$ be
the perturbed points of Example 1.1. Since the preprocessing phase
does not change the perturbed points, we have $\overline{\cal
P}=\widehat{\cal P}$, corresponding to the exact normal set
$\overline{\cal N}_0 =\{1,\; y,\; x\}$. Since the relative
residuals of the terms  $y$ and $x$ computed on $\overline{\cal
P}$ are equal to  $0.37$ and $0.0068$, respectively, we choose
$\epsilon > 0.0068$, to obtain a set $\overline{GB}$ different
from the exact G-basis of ${\cal I}(\overline{\cal P})$.

Choosing  $0.0068 <\epsilon < 0.37$, the polynomial
$x-2.05y+1.0{\overline 6}$ is inserted into $\overline{GB}$ and
the term  $y^2$ is analyzed. Since the relative residual of $y^2$
is equal to $0.0824$ we have two possibilities.

If $0.0824 <\epsilon < 0.37$, the nonzero relative residual of
$y^2$ is treated as if it were equal to zero and the following set
of polynomials is computed:
\begin{eqnarray*}
{\overline {GB}_1=} \; & & \left\{
\begin{array}{l}
x-2.05y+1.0\overline 6 \\
y^2-4y+3.\overline 3.
\end{array}
\right.\end{eqnarray*}
 Since $\overline{GB}_1$  vanishes at $\overline{\cal P}^+_1=\{(4.7072,
\;2.8165),\;(1.3595,\;1.1835)\}$ which differs in cardinality from
$\overline{\cal P}$, we conclude that $\overline {GB}_1$ is not a
good approximation to $GB$ and that the threshold has not been
suitably chosen.

If $0.0068 <\epsilon <0.0824 $, then the polynomial
$y^2-4y+3.\overline 3$ does not belong to $\overline{GB}_2$  and
the term $y^2$ is a new element of the normal set. So in this case
the ABM algorithm computes the set $\overline{GB}_2$:
\begin{eqnarray*}
\overline{GB}_2 =  & & \left\{
\begin{array}{l}
x-2.05y+1.0{\overline 6} \\
y^3-6y^2+11y-6
\end{array}
\right.\end{eqnarray*}
G-basis of the ideal of points $\overline{\cal P}_2^+=\{(0.98{\overline
3},\; 1),\;(3.0{\overline 3},\; 2),\; (5.08{\overline 3},\;3)\}$,
very near to the set of points $\overline{\cal P}$, and so we
consider the threshold suitably chosen. Note that
$\overline{GB}_2$ is very similar to the G-basis of the ideal of
the exact points of Example 1.1.\hfill$\diamondsuit$

\vspace{3 mm} {\bf Remark.} If the ABM algorithm  cannot compute
normal sets which satisfy  properties P1--P3, then  we are dealing
with ``degenerate" cases. The numerical tests suggest that, in the
degenerate examples, the relative least square residuals for the
exact normal sets associated to all the admissible input points
are small. Unfortunately, since it is rather difficult, if not
impossible, to distinguish the case where the small residuals
computed on the preprocessed points are caused by data
perturbations from the case where they are related to degenerate
problems, we cannot check, {\it a priori}, if we are dealing with
degenerate cases.

\section{ Comparison with a reference basis }

In this section we estimate the difference between a reference
basis and the set $\overline{GB}$ computed, given the term
ordering $\sigma$, by the ABM algorithm working on the
preprocessed points $\overline{\cal P}$. We assume that the value
of the threshold $\epsilon$ has been suitably chosen  and that
``non-degenerate" problems are analyzed, that is the computed set
$\overline{\cal N}$ is a reference normal set and there exists a
reference $\sigma$-Gr\"obner basis $GB$ with normal set
$\overline{\cal N}$. We show that $\overline{GB}$ is a good
approximation to a reference basis $GB$.

First of all, note that if at each step of the ABM algorithm the
analyzed relative residual is equal to zero or greater than the
threshold, then the ABM algorithm works on the preprocessed points
$\overline{\cal P}$ like the BM algorithm. In this case the set
$\overline{GB}$ is the exact G-basis of the ideal ${\cal
I}(\overline{\cal P})$ and so, since  $\overline{\cal P}$ is an
admissible input set, $\overline{GB}$ coincides with a reference
basis.

In the general case, it can happen that, at some step of the ABM
algorithm, the relative residual associated to $\overline t$ is
greater than zero and less than or equal to $ \epsilon$. In this
case the polynomial $\overline g = \overline t - \sum_{i=1}^k
\overline c_i  \overline t_i$ is the new element of the set
$\overline{GB}$, even though it does not vanish at the points
$\overline{\cal P}$, since
 $$\overline g(\overline{\cal P}) =  \overline t(\overline{\cal P})
- \sum_{i=1}^k \overline c_i \overline t_i(\overline{\cal
P})=\overline  r,   $$%
where $\overline r$ is the nonzero residual of (\ref{pbminq}), and
so $\overline{GB}$ is not the G-basis of the ideal ${\cal
I}(\overline{\cal P})$. Nevertheless,  since  $\overline{\cal N}$
is also the normal set of  a reference basis $GB$, $\overline{GB}$
structurally corresponds to $GB$. In addition, in the following we
prove that the polynomials of  $\overline{GB}$ and $GB$ have
similar coefficients. For these reasons, as mentioned earlier, the
computed set $\overline {GB}$ can be considered a good
approximation to a reference basis.

Let $\cal P$ be the admissible input points which are zeros of the
polynomials of a reference basis $GB$. By construction, if
$\overline g \in \overline{GB}$ corresponds to $g\in GB$, the
vectors $\overline c$ and  $c$ of the coefficients of $\overline
g$ and $g$  satisfy respectively
\begin{eqnarray*}
\overline M_k(\overline{\cal P}) \overline c  = \overline
t(\overline{\cal P})- \overline r \;\;{\rm and} \;\; \overline M_k({\cal P})
c = \overline t({\cal P}),
\end{eqnarray*}
where $\overline r$ is the residual of the least square problem
(\ref{pbminq}) using the preprocessed points. From the analysis of
the sensitivity of a least square problem (\cite{DBA},
\cite{GVL}), since the residual of the system
(\ref{sistema_esatto}) is zero, we have
\begin{eqnarray}\label{upbound1}
\frac{\|\overline c - c\|_2}{\|c\|_2} \le \frac{K_2(\overline
M_k(\overline{\cal P}))} {1-\epsilon_M K_2(\overline
M_k(\overline{\cal P}))}[\epsilon_t+\epsilon_M],
\end{eqnarray}
where
\begin{eqnarray*}
\epsilon_t=\frac{\|\overline t(\overline {\cal
P})-\overline t({\cal P})\|_2}{\|\overline t({\cal P})\|_2},\;\;
\epsilon_M=\frac{\|\overline M_k(\overline{\cal P})-\overline
M_k({\cal P})\|_2}{\|\overline M_k({\cal P})\|_2},
\end{eqnarray*}
and $K_2(\overline M_k(\overline{\cal P}))$ is the 2-condition
number of the matrix $\overline M_k(\overline{\cal P})$
\cite{BCM}. It is well known \cite{BCM} that, in general, the
upper bound (\ref{upbound1}) overestimates the sensitivity of the
linear systems. Nevertheless,  the formula (\ref{upbound1})
suggests that $\epsilon_t$ and $\epsilon_M$ can be interpreted as
a measure of goodness of approximation. If $\delta_R$ is the
maximum relative error which perturbs the coordinates of the
points $\overline{\cal P}$ with respect to $\cal P$, we show, in
Appendix B, that
$$\epsilon_t \le Deg(\overline t) \delta_R \;\;\;\; {\rm and}\;\;\;\; \epsilon_M
 \le \sqrt{k}\,d_M \delta_R,$$
where $Deg(\overline t)$ is the degree of the term $\overline t$
and $d_M$ is the maximum degree of the terms $\overline
t_1,\dots,\overline t_k$. Since both sets $\overline{\cal P}$ and
$\cal P$ are admissible data perturbations, the differences
between their coordinates is not large and so the value $\delta_R$
is small. We can conclude that, with a suitable choice of the
threshold $\epsilon$, the algorithm ABM, working on the
preprocessed  points $\overline{\cal P}$,  computes, for
non-degenerate problems, a good approximation to a reference
basis.

\vspace{3 mm} {\bf Remark.} Note that, under small perturbations
of $\overline{\cal P}$, the matrix $\overline M_k(\overline{\cal
P})$ of the system (\ref{pbminq}) preserves its full rank, since
$\overline{\cal N}$ is a reference normal set, but its elements
slightly change. For this reason,  choosing different admissible
input sets, we obtain sets of polynomials which differ only in
their coefficients, in a continuous way, and  which do not present
structural differences.

\section{Extended basis and pseudozeros}

In his book ``Numerical Polynomial Algebra" \cite{Ste04} and in
several papers (\cite{Ste97}, \cite{Ste199}, \cite{Ste299},
\cite{Ste00}, \cite{Ste01}) Stetter analyzes some aspects of
Polynomial Algebra from the numerical
point of view, when perturbed data and floating point arithmetic are
used. In this
section we compare the set $\overline {GB}$ with the {\it extended
basis} developed by Stetter. Moreover, we show that, in the
terminology of Stetter, $\overline{\cal P}$ is a pseudozero set
for the polynomials of $\overline{GB}$.

\subsection{The extended basis}
In \cite{Ste97}, Stetter analyzes the problem of computing the
G-basis of an ideal determined by a system of polynomials with
perturbed coefficients, pointing out that the G-basis can be
structurally altered by the data uncertainty. The numerical
instabilities occur, analogously to the case of ideals of points,
when the zeros of the perturbed polynomials are ``near" to a set
of points associated to a normal set ${\cal N}^\rho$ which differs
from the perturbed normal set $\widehat{\cal N}^\rho$. In order to
solve the numerical instabilities, Stetter introduces the notion
of {\bf extended basis} $GB_E$, that is an approximation to the
G-basis, computed using the normal set ${\cal N}^ \rho$ instead of
$\widehat{\cal N}^\rho$. The computed set $GB_E$ is not a G-basis
in general, but, if there exists a G-basis $GB$ structurally
corresponding to $GB_E$  such that for each polynomial $g_E \in
GB_E$ there is a polynomial $g \in GB$ with ``similar"
coefficients, then $GB_E$ is a ``numerically stable" basis of the
ideal associated to the input polynomials (see theorem 4.1 in
\cite{Ste97}).

Modifying the technique proposed by Stetter, it is possible to
compute the extended basis for ideals of perturbed points too.
Given the set of $m$ points $\overline{\cal P}$, the normal set
${\cal N}^\rho$ coincides with the normal set $\overline{\cal
N}=\{\overline t_1,\dots,\overline t_m\}$ computed by the ABM
algorithm, because $\overline{\cal N}$ characterizes all the
admissible input sets. Following Stetter's method, each polynomial
$g_E\in GB_E$ corresponds to a power product $\overline t \in
{\cal C}[\overline{\cal N}]$ and can be written as
\begin{eqnarray}\label{esteso}
g_E=\overline t - \sum_{j=1}^m c^E_j\overline t_j,
\end{eqnarray}
using all the terms of $\overline{\cal N}$. The vector $c_E$ of
the  coefficients of $g_E$ is the solution of the system
\begin{eqnarray}\label{mat_estesa}
\overline M_m(\overline{\cal P}) c_E = \overline t(\overline{\cal
P}), \;\;\;{\rm with}\;\;\; \overline M_m(\overline {\cal P})
=\left[ \overline t_1(\overline{\cal P}), \dots, \overline
t_m(\overline{\cal P}) \right ],\end{eqnarray}
 and so $g_E$ vanishes at $\overline{\cal P}$.

 Since for each
$\overline t \in {\cal C}[\overline{\cal N}]$ there exist both  a
polynomial $\overline g \in \overline{GB}$ and a polynomial $g_E
\in GB_E$, there exists an one-to-one correspondence between
$\overline{GB}$ and $GB_E$. In addition, since the vector
$\overline c$ of the coefficients of $\overline g$ is the least
square solution of the problem (\ref{pbminq}), we have that
$$ \overline g(\overline{\cal P}) =\overline t(\overline{\cal P})
- \overline M_k(\overline{\cal P}) \overline c, \;\;\;{\rm and \;
from\; (\ref{mat_estesa})},\;\;\;  \overline g(\overline{\cal P})
= \overline M_m(\overline{\cal P}) c_E -\overline
M_k(\overline{\cal P}) \overline c. $$
 Since the matrix
$\overline M_k(\overline{\cal P})$ consists of the first $k$
columns of $\overline M_m(\overline{\cal P})$,   we have
$$ \overline g(\overline{\cal P}) = \overline M_m(\overline{\cal P}) c_E -
\overline M_m(\overline{\cal P})\left[
\begin{array}{c} \overline c \\ 0\end{array}\right ] = \overline M_m(\overline{\cal
P}) \Delta c,\;\;{\rm where}\;\; \Delta c = c_E-
\left[\begin{array}{c}\overline c \\ 0\end{array}\right ].$$ Since
$\Delta c=(\overline M_m(\overline{\cal P}))^{-1} \overline
g(\overline{\cal P})$ and $c_E=(\overline M_m(\overline{\cal
P}))^{-1} \overline g(\overline{\cal P})$, we obtain \cite{GVL}
\begin{eqnarray*}
\frac{\|\overline g(\overline{\cal P}) \|_2}{\sigma_{\max}} \le
\|\Delta c\|_2 \le \frac{\|\overline g(\overline{\cal
P})\|_2}{\sigma_{\min}}\;\;\;{\rm and }\;\;\; \frac{\|\overline t(\overline{\cal
P})\|_2}{\sigma_{\max}} \le \|c_E\|_2 \le \frac{\|\overline
t(\overline{\cal P}) \|_2} {\sigma_{\min}},
\end{eqnarray*}
where $\sigma_{\max}$ and $\sigma_{\min}$ are respectively the
minimum  and the maximum singular values of $\overline
M_m(\overline{\cal P})$. Since the ABM algorithm computes an
element of $\overline{GB}$ if the relative residual of the least
square problem (\ref{pbminq}) is less than the threshold
$\epsilon$, we can derive the following  upper bound for the
differences of the coefficients of the elements of $\overline
{GB}$ and of $GB_E$:
$$\frac{\|\Delta c \|_2}{\|c_E \|_2} \le \frac{\|\overline g(\overline{\cal P})\|_2}{\|
\overline t(\overline{\cal
P})\|_2}\frac{\sigma_{\max}}{\sigma_{\min}} \le \epsilon
K_2(\overline M_m(\overline{\cal P})),$$ where $K_2(\overline
M_m(\overline{\cal P}))=\sigma_{\max}/\sigma_{\min}$ is the
2-condition number of $\overline M_m(\overline{\cal P})$. Since in
general the previous upper bound overestimates the relative error
on the coefficients, we can conclude that $\overline{GB}$ is
``near" to the extended basis $GB_E$.

\vspace{3 mm} {\bf Remark.} The power product $\overline t \in
{\cal C}[\overline{\cal N}]$ is probably not the leading term of the
polynomial $g_E$ defined by (\ref{esteso}). In fact, if $\overline
t_k <  \overline t < \overline t_{k+1}$, since in general
$\overline t(\overline{\cal P})$ is not in the span of the first $k$
columns of $\overline M_m(\overline{\cal P})$ in the classical
sense, then the last $m-k$ components of $c_E$ are different from
zero. In contrast, by construction, $\overline t$
is the leading term of $\overline g \in \overline{GB}$ and so,
with respect to the extended basis, $\overline{GB}$ has the
advantage that $LT(\overline{GB})={\cal C}[\overline{\cal N}]$,
that is $LT(\overline{GB})$ coincides with the leading term of a
reference basis.

\subsection{Pseudozeros} Even if sometimes $\overline{\cal
P}$ is a zero set for the polynomials of $\overline{GB}$, in
general the elements of $\overline{GB}$ do not vanish at
$\overline{\cal P}$, but we show that the preprocessed points are
pseudozeros for the polynomials of $\overline{GB}$.

Intuitively, given a tolerance $\delta$, a point $p$ is a
pseudozero of a polynomial $g$ if it is the exact zero of a
polynomial $\hat g$ whose coefficients differ from those of $g$ by
less than $\delta$. This idea, formalized by Stetter in
\cite{Ste199}, can be generalized for a set of points.
\begin{defin}
Let $S$ be  a set of power products.

\noindent A point $p$ is a pseudozero of $g=\sum_{x^j\in S} a_j
x^j $ with respect to $S$ and  tolerance $\delta$ if $p$ is an
exact zero of  some $\tilde g \in N_S(g, \delta)$, where
$$N_S(g,\delta)=\{ \tilde g : \tilde g (x)= \sum_{x^j\in S} \tilde a_j x^j,
|\tilde a _j - a_j | < \delta \}.$$
 A set of point $ \cal P$ is a set of
pseudozeros of $g$ with respect to $S$ and tolerance $\delta$ if
it is a set of exact zeros of  some $\tilde g \in N_S(g,\delta)$.

\noindent
 A set of point $ \cal P$ is a set of
pseudozeros  for a system of $k$ polynomials $g_1,\dots,g_k$,
whose supports belong to $S$, with respect to $S$ and tolerances
${\cal D}=\{\delta_1,\dots,\delta_k\}$ if it is a set of
pseudozeros with tolerance $\delta_i$ for each polynomial $g_i$,
$i=1\dots k$, of the system.
\end{defin}
The set $\overline{\cal P}$ of the preprocessed points is a
pseudozero set for $\overline{GB}$ with respect to
$S=\overline{\cal N} \cup {\cal C[\overline N]}$, as shown in the
following theorem.
\begin{theo}
The points of  $\overline{\cal P}$ are pseudozeros of the system
(\ref{sistema}) consisting of the polynomials
$\overline{GB}=\{\overline g_1,\dots,\overline g_m\}$ with respect
to $S$ and with tolerances ${\cal D} =\{ \delta_1, \dots,
\delta_m\}$, where
$$S=\overline{\cal N} \cup {\cal C[\overline N]}\;\;and\;\;\;\delta_i=\frac{\|\overline g_i({\overline{\cal P}})\|_2}{\sigma_{\min}}.$$
\end{theo}

\noindent{\bf Proof.} Given a polynomial $\overline g_i =
\overline t - \sum_{j=1}^k \overline c_j \overline t_j\;\in
\overline{GB}$, there exists $g_E = \overline t - \sum_{j =1}^m
c^E_j \overline t_j\;\in{GB_E}$, vanishing at $\overline{\cal P}$,
such that
\begin{eqnarray*}
|c^E_j - \overline c_j| \le    \|\Delta c\|_2 \le \frac{
\|\overline g_i(\overline {\cal P})\|_2}{\sigma_{\min}}=\delta_i.
\end{eqnarray*}
Since the supports of $\overline g$ and $g_E$ are contained in
$S=\overline{\cal N} \cup \,{\cal C[\overline N]}$ and $g_E$
belongs to $N_S(\overline g_i, \delta_i)$ we have that
$\overline{\cal P}$ is a pseudozero set of $\overline{GB}$.
\hfill{$\diamondsuit$}

Each value $\delta_i$ is small with respect to the norm of the
vector $\overline c$ of the coefficients of $\overline g_i$. In
fact, if $\overline g_i \in \overline{GB}$ and $g_E \in GB_E$ have
similar coefficients,  from the results of the previous section we
have
$$ \delta_i \le \epsilon K(\overline M_m(\overline{\cal
P}))\|c_E\| \;\;\; {\rm and \; so} \;\;\; \delta_i \le \epsilon
K(\overline M_m(\overline{\cal P}))\|\overline c\|+O(\epsilon^2).
$$
The conclusion follows, since the previous upper bound is, in
general, an  overestimation and $\epsilon$ is a small value of the
threshold.

\section{Numerical examples}

In this section we present two numerical examples which illustrate
how the ABM algorithm approximates a reference basis. From the
perturbed points, after the preprocessing phase, the ABM algorithm
computes in Example 6.1 a reference basis and in Example 6.2 an
approximation to it. In these examples we use the term ordering
DegLex; in addition, the coordinates of the points and the
coefficients of the polynomials  are displayed with a finite
number of digits, but all computations are performed in exact
arithmetic using CoCoA 4.2 \cite{CoCoA}. Moreover, in the
following examples the threshold in the test of the ABM algorithm
has been suitably chosen. A more exhaustive survey of numerical
experiments will be presented in a forthcoming paper.

\vspace{3 mm} \noindent{\bf Example 6.1}

\noindent
Given the perturbed points $\widehat{\cal P}$ known with
uncertainty $s_0=0.1$,
\begin{eqnarray*}
\widehat{\cal P}=\{(-2.45,-3.6 ), \;(-0.53,-1.45), \;(1.5, 0.45 ),
\;(3.5, 2.5)\},
\end{eqnarray*}
 the exact DegLex-Gr\"obner basis  of the ideal ${\cal I}(\widehat{\cal P})$
is
\begin{eqnarray*}
\widehat{GB}=\left\{\begin{array}{l}
xy-0.97550y^2+1.45450x-2.48031y-1.54307\\
x^2-0.95161y^2+0.89208x-2.94110y-2.07192\\
y^3+0.64549y^2+91.31103x-98.56564y-92.83385
\end{array}\right.
\end{eqnarray*}

The exact G-basis $\widehat{GB}$ cannot be a reference basis, since there is a small relative residual associated to a power product of its normal set $\widehat{\cal N}$, as pointed out in the following table.
\begin{table}[h]
\begin{center}
\begin{tabular}{c|cccc}
$\widehat{\cal N} $& $1$ & $y$& $x$ & $y^2$ \\
\hline
Rel. Res.  &  & $0.97$ & ${\bf 0.03}$ &$0.45$
\end{tabular}
\end{center}
\end{table}

The small relative residual associated to the power product $x$
means that the perturbed points are ``almost aligned" and so a
normal set corresponding to the ``configuration" of aligned points
 turns out to be more meaningful and stable than $\widehat{\cal N}$. The ABM
algorithm can detect this situation. In fact, using for example the threshold
$\epsilon=0.1$ and working on the preprocessed points
$\overline{\cal P}$
\begin{eqnarray*}
\overline{\cal P}=\{(-2.475,-3.55 ), \;(-0.49,-1.475),
\;(1.475, 0.49), \;(3.55, 2.475)\},
\end{eqnarray*}
the ABM algorithm computes the  set $\overline {\cal
N}=\{1,y,y^2,y^3\}$ and the set of polynomials
\begin{eqnarray*}
\overline{GB}=\left\{\begin{array}{l}
x-0.99979y-1.02989\\
y^4+2.06y^3-8.45012y^2-9.43141y+6.35026.
\end{array}\right.
\end{eqnarray*}
$\overline{\cal N}$ is a reference normal set, since the relative
residuals associated to its terms are  greater than $0.22$ and the
relative residuals associated to $x$ and $y^4$ are less than $
0.02$. In addition, the set $\overline{GB}$ is a reference basis,
since it is the DegLex-Gr\"obner basis of the ideal ${\cal I}(\overline{\cal
P}^+)$, with normal set
$\overline{\cal N}$, where $\overline{\cal P}^+$ is the following admissible
input set computed by Newton iteration
\begin{eqnarray*}
\overline{\cal P}^+=\{(-2.519,-3.55),\,(-0.444,-1.475),\,(
1.519 ,0.49),\,( 3.504, 2.475)\}.\;\;\;\; \diamondsuit
\end{eqnarray*}

\vspace{3 mm}\noindent{\bf Example 6.2}

\noindent Let $\widehat{\cal P}$ be the set of the perturbed
points with uncertainty $s_0=0.1$,
\begin{eqnarray*}
\widehat{\cal P}=& & \{(0, 4.1),\;(0.05,-4.03 ),\;(3.1, 0.1),\;
(-3, 0.03),\\ &&(2.37,2.5),\;(-2.4,-2.33),\;
(2.31,-2.486),\;(-2.4,2.4)\}.
\end{eqnarray*}
 The exact DegLex-Gr\"obner basis of the ideal ${\cal I}(\widehat{\cal P})$ is
\begin{eqnarray*}
\widehat{GB}=\left\{\begin{array}{ll}%
x^2y &-0.00164xy^2+0.52911y^3+0.15329x^2-0.01032xy \\
& +0.05231y^2-0.02767x-8.75015y-1.47114, \\
x^3 & +0.60405xy^2+0.00120y^3-1.38130x^2-0.04581xy \\
& -0.72231y^2-9.17153x+0.03536y+11.91429,\\
y^4&-2.50589xy^2-0.34459y^3-178.88212x^2+1.01537xy\\
&-117.12255y^2+17.87987x+11.35992y+1663.42803,\\
xy^3&-0.09269xy^2+0.04595y^3+1.13505x^2-5.90939xy\\
&+0.65885y^2+0.28526x-1.06888y-9.86020
\end{array}\right.
\end{eqnarray*}

The exact normal set $\widehat{\cal N}$ of $\widehat{GB}$ cannot
be a reference normal set, since the relative residual associated
to its power product $x^2$ is very small, as reported in the
following table.

\begin{table}[h]
\begin{center}
\begin{tabular}{c|ccccccccc}
$\widehat{\cal N} $& $1$ & $y$& $x$ & $y^2$ &$xy$ & $x^2$&$y^3$&$xy^2$ \\
\hline
Rel. Res.  &  & $0.99$ & $1$ &$0.64$& $0.99$& ${\bf 0.03}$ & $ 0.39$ & $0.57$ \\
\end{tabular}
\end{center}
\end{table}

\noindent Since $\widehat{\cal N}$ can change structurally for
small data perturbations,  $\widehat{GB}$ is not a reference
basis. The ABM algorithm, using the threshold $\epsilon=0.1$ and
working on the preprocessed points $\overline{\cal P}$
\begin{eqnarray*}
\overline{\cal P}=&&\{(0,4.065),\;(0,-4.065),\;(3.05,0),\;(-3.05,0),\\
&&(2.405 ,2.405),\;(-2.405 ,-2.405),\;(2.405, -2.405 ),\;( -2.405
, 2.405)\}
\end{eqnarray*}
computes a set of power products $\overline{\cal N}=\{1,y,x,y^2,xy,y^3,xy^2,y^4\}$ which is a normal set associated to
relative residuals greater than $0.22$. Since all the relative
residuals are sufficiently greater or less than the threshold, $\overline{\cal
N}$ turns out to be a reference normal set. Moreover, the ABM
algorithm computes the set of polynomials
\begin{eqnarray*}
\overline{GB}=\left\{\begin{array}{l}
x^2+0.55840y^2-9.13935, \\
xy^3-5.78402xy,\\
y^5-22.30825y^3+95.57653y,
\end{array}\right.
\end{eqnarray*}
whose zero set $\{(3.02313,0),(-3.02313,0)\}$ is not an admissible
input set. Nevertheless, Newton method applied to the function
$\overline F$, defined by (\ref{funzione}), follows  descent
directions for $\|\overline F\|$. For this reason the Newton
iteration computes the admissible input set
\begin{eqnarray*}
\overline{\cal P}^+&=&\{(0,4.065),\;(0,-4.065),\;(3.02,0),\;(-3.02,0),\\
& &(2.41 ,2.41),\;(-2.41 ,-2.41),\;(2.41, -2.41 ),\;( -2.41,
2.41)\}.
\end{eqnarray*}
 This fact suggests that the threshold has been suitably chosen and
that $\overline{\cal N}$ and $\overline{GB}$ are a good
approximation to the exact case. In fact, the set $\overline{GB}$
approximates the set $GB$
\begin{eqnarray*}
GB=\left\{\begin{array}{l}
x^2+0.5625y^2-9,\\
xy^3-5.76xy,   \\
y^5-21.76y^3+92.16y
\end{array}\right.
\end{eqnarray*}
which is a reference basis, since it is a DegLex-Gr\"obner basis,
with normal set $\overline{\cal N}$, of the ideal of the
admissible input points
\begin{eqnarray*}
\{(0,4),(0,-4),(3,0),(-3,0),(2.4,2.4),(-2.4,-2.4),(2.4,-2.4),(-2.4,2.4)\}.&&\\
\diamondsuit &&
\end{eqnarray*}

\section*{Appendix A: the preprocessing phase}
In this appendix we summarize the preprocessing  phase of the
perturbed points $\widehat{\cal P}$. This procedure finds the
coordinates which differ by less than the error estimation and
replaces them by a single representative value. For this reason
the set $\overline{\cal P}$ of the preprocessed points is
equivalent to $\widehat {\cal P}$ from the computational point of
view, even if it can happen that the cardinality of $\overline
{\cal P}$ is less than the cardinality of $\widehat {\cal P}$. In
order to preserve or to get back the ``geometrical symmetries"  of
the exact points, the preprocessing algorithm works on the
absolute values of the coordinates of all the points, assembled
together in a set $Y$, and it computes a preprocessed set
$\overline Y$. The preprocessed points $\overline{\cal P}$ are
built from the set $\overline Y$ in an obvious way.

Let $s_0$ be an estimate of the errors on the data and let
$Y=\{y_1,\dots,y_n\}$ be  the set of the absolute values of the
coordinates of the input points, sorted in non-decreasing order.
If $(k+1)$ elements $y_{j},\dots,y_{j+k}$ of $Y$ are such that the
intersection
$$X^k_j=\bigcap_{i=j}^{j+k}[y_i-s_0,\;y_i+s_0]=[y_{j+k}-s_0,y_j+s_0],$$
is not empty, then any element $\overline y_{j} \in X^k_j$ can
represent, from the computational point of view, the values
$y_{j},\dots,y_{j+k}$, since it differs from them by less than
$s_0$. For this reason we replace $y_j,\dots,y_{j+k}$ with a value
$\overline y_j\in X_j^k$ in the set $\overline Y$ of the
preprocessed coordinates.

The preprocessing algorithm for building the set $\overline Y$ can
be described as follows. First of all, each element of $Y$ less
than $s_0$ is removed from $Y$ and the value $0$ is inserted into
$\overline Y$. After this step, the algorithm computes for each
element $y_j \in Y$ the largest non empty intersection $X_{j}^k$,
and then it processes the set $X_j^k$ which contains the maximum
number of elements of Y, that is the set with the maximum index
$k$. If there are several intersections with this property, the
set with the minimum index $j$ is chosen. If $X_j^k$ is the
intersection to process, all the elements $y_{j},\dots,y_{j+k}\in
X_j^k$ are removed from $Y$ and the middle point $\overline y_j$
of $X_j^k$ is inserted in $\overline Y$.

The preprocessing algorithm ends when $Y$ is empty or when it
contains elements which differ by more than $s_0$. In this case
such values do not require any preprocessing treatment and so they
are removed from $Y$ and inserted directly into $\overline Y$.

Even though  it can happen that $[\overline y_i-s_0,\; \overline
y_i+s_0] \cap [\overline y_{i+1}-s_0,\; \overline y_{i+1}+s_0]\neq
\emptyset$, we decided not to repeat the preprocessing phase on
the set $\overline Y$, since  all the new  coordinates  differ by
more than $s_0$. In fact, if the set $X_{j}^k$ is processed,
\begin{eqnarray*}
[y_{j-1}    -s_0,\; y_{j-1}+s_0] \cap X_{j}^k=\emptyset, \;\;\;\;
[y_{j+k+1}  -s_0,\; y_{j+k+1}+s_0] \cap X_{j}^k =\emptyset,
\end{eqnarray*}
and so the middle point $\overline y_{j}$ of  $X_{j}^k$ differs
from $y_{j-1}$ and $y_{j+k+1}$ by more than $s_0$. Since
$y_{j+k+1}-y_{j-1} > 2s_0$, the preprocessing phase works
separately on the sets $Y_1=\{y_1,\dots,y_{j-1}\}$ and $
Y_2=\{y_{j+k+1},\dots,y_n\}$. The values obtained by preprocessing
$Y_1$ and $Y_2$ are, respectively, less than the maximum of $Y_1$
and greater than the minimum of $Y_2$, and so they differ from
$\overline y_j$ by more than  $s_0$. Analogously, analyzing the
sets $Y_1$ and $Y_2$ separately, we can show that all the
preprocessed values differ to each other by more than $s_0$, that
is they are ``well separated".
This result points out that it is computationally better
to use the set $\overline{\cal P}$ than the set $\widehat {\cal
P}$ as input points.

\section*{Appendix B: an upper bound for $\epsilon_t$ and $\epsilon_M$}
In this section we present an upper bound for the values
$\epsilon_t$ and $\epsilon_M$,  useful for estimating how well the
set $\overline{GB}$ approximates a reference basis.

First of all, we analyze the sensitivity of a term evaluated at a
perturbed point. Let $t=x_1^{j_1}\dots x_s^{j_s}$ be a power
product belonging to $\Reali[x_1\dots x_s]$ and let
$Deg(t)=j_1+\dots+j_s$ be its degree. If $\overline p =(
\overline p_1,\dots,\overline p_s)$ is a perturbation of the point
$p =(  p_1,\dots, p_s)$ such that $\overline p_i=p_i(1+\delta_i),
\; |\delta_i| \le \delta_p,\; i=1\dots s, $ then
\begin{eqnarray*}
t(\overline p)= \overline p_1^{j_1} \dots  \overline p_s^{j_s}=p_1^{j_1}(1+\delta_1)^{j_1} \dots  p_s^{j_s}(1+\delta_s)^{j_s}= t(p) (1+\delta_1)^{j_1}\dots (1+\delta_s)^{j_s}.
\end{eqnarray*}
By a  first-order error analysis, ignoring the errors of higher order, we obtain
\begin{eqnarray*}
|t(\overline p)-t(p)| \le |t(p)| (j_1 |\delta_1| + \dots j_s |\delta_s|) \le
 |t(p)| Deg(t) \delta_p.
\end{eqnarray*}
Now, we can easily upper bound  $\epsilon_t$ as follows. Let
$\overline{\cal P}=\{\overline p_1,\dots,\overline p_m\}$ be a set
of points perturbation of ${\cal P}=\{p_1,\dots,p_m\}$, such that
$p_i=(p_1^{(i)},\dots,p_s^{(i)})$ and $\overline
p_i=(p_1^{(i)}(1+\delta_1^{(i)}),\dots,p_s^{(i)}(1+\delta_s^{(i)}))$,
where $|\delta_j^{(i)}| \le \delta_R $, $i=1\dots m$ and $j=1\dots
s$. Then we have
\begin{eqnarray*}
\|t(\overline{\cal P})-t({\cal P})\|_2^2  \le \sum_{i=1}^m
t(p_i)^2 Deg(t)^2 \delta_R^2 \le Deg(t)^2 \delta_R^2 \|t({\cal P})
\|_2^2 \;\;
 {\rm that \; is} \;\; \epsilon_t \le  Deg(t) \delta_R.\end{eqnarray*}
Analogously, we can also upper bound $\epsilon_M$. Let $M({\cal
P})$ and $M(\overline{\cal P})$ be the matrices whose columns are
generated by the terms $t_1 \dots t_k$ evaluated, respectively, at
$\cal P$ and $\overline{\cal P}$. If $d_M= \max \{Deg(t_1),
\dots,Deg(t_k)\}$ is the maximum degree of the terms $t_1 \dots
t_k$, we obtain
$$\|M(\overline{\cal P})- M({\cal P})\|_F ^2 = \sum_{j=1}^k \|t_j(\overline{\cal P})
 - t_j({\cal P})\|_2^2 \le \delta_R^2 d_M^2 \sum_{j=1}^k \|t_j({\cal P})\|_2^2 = d_M^2 \delta_R^2 \|M({\cal P})\|_F^2,$$
where $\| \phantom{0}\|_F$ is the Frobenius matrix norm. Since, given an $m\times k$ matrix $A$, we have  $\|A\|_2 \le \|A\|_F \le \sqrt{k}\|A\|_2$, then
\begin{eqnarray*}
\epsilon_M = \frac{\|M(\overline{\cal P}) - M({\cal P})\|_2}{\|M({\cal P})\|_2}\le
\sqrt{k}\frac{\|M(\overline{\cal P}) - M({\cal P})\|_F}{\|M({\cal P})\|_F} \le \sqrt{k} \, d_M \delta_R.
 \end{eqnarray*}

\vspace{3 mm}\noindent {\bf Acknowledgements.} The author would
like to thank John Abbott for his useful and constructive remarks.

\end{document}